\documentclass{article}

\usepackage{amssymb, amsmath, latexsym, mathrsfs, graphicx, stmaryrd, amsthm, psfrag}
\usepackage{url}
\newtheorem{thm}{Theorem}[section]
\newtheorem{lemma}[thm]{Lemma}
\newtheorem{cor}[thm]{Corollary}

\newtheorem{prop}[thm]{Proposition}

\theoremstyle{definition}
\newtheorem{exm}[thm]{Example}
\newtheorem{definition}[thm]{Definition}

\newcommand{\N}{\mathbb{N}}
\newcommand{\Z}{\mathbb{Z}}

\newcommand{\R}{\mathbb{R}}
\newcommand{\C}{\mathbb{C}}

\makeatletter
\def\revddots{\mathinner{\mkern1mu\raise\p@
\vbox{\kern7\p@\hbox{.}}\mkern2mu
\raise4\p@\hbox{.}\mkern2mu\raise7\p@\hbox{.}\mkern1mu}}
\makeatother
\def \vec#1{\mathbf{#1}}
\def \abs#1{\lvert #1 \rvert}
\renewcommand{\mod}[1]{(\text{mod }#1)}

\begin{document}

\title{Vertices of Gelfand-Tsetlin Polytopes\footnote{ Research
supported by NSF Grant DMS-0309694 and by NSF VIGRE Grant No. 
DMS-0135345.}}

\author{ Jes\'us A. De Loera and Tyrrell B. McAllister}

\maketitle

\noindent \emph{This paper is dedicated to Louis Billera on the
occasion of his sixtieth birthday.}

\vskip .5cm

\noindent\textbf{Abstract:} This paper is a study of the polyhedral geometry of Gelfand-Tsetlin
patterns arising in the representation theory $\mathfrak{gl}_n \C$ and algebraic combinatorics.
We present a combinatorial characterization of the vertices and a method
to calculate the dimension of the lowest-dimensional face containing a given
Gelfand-Tsetlin pattern.

As an application, we disprove a conjecture of Berenstein and Kirillov
\cite{berensteinkirillov} about the integrality of all vertices of the Gelfand-Tsetlin polytopes. We can
construct for each $n\geq5$ a counterexample, with arbitrarily increasing denominators
as $n$ grows, of a non-integral vertex. This is the first infinite family of
non-integral polyhedra for which the Ehrhart counting function is still a polynomial.
We also derive a bound on the denominators for the non-integral vertices
when $n$ is fixed.

\section{Introduction}

Many authors have recently observed that polyhedral geometry plays a special
role in combinatorial representation theory (see for example
\cite{berensteinzelevinsky,kingetal, kirillov2,knutsontao, nakashimazelevinsky},
and the references within). In this note we study the polyhedral
geometry of the so-called Gelfand-Tsetlin patterns, which arise
in the representation theory of $\mathfrak{gl}_n \C$ and the study of
Kostka numbers.

For each $n \in \N$, let $X_n$ be the set of all triangular arrays
$(x_{ij})_{1\le i\le j\le n}$ with $x_{ij} \in \R$.  Then $X_n$
inherits a vector space structure under the obvious isomorphism $X_n
\cong \R^{\frac{n(n+1)}{2}}$.

\begin{definition}{\label{GT-patterns}}
A \emph{Gelfand--Tsetlin pattern} or \emph{GT-pattern} is a triangular
array $(x_{ij})_{1\le i\le j\le n} \in X_n$ satisfying the inequalities
\begin{itemize}
\item $x_{ij} \geq 0$, for $1 \leq i \leq j \leq n$; and
\item $x_{i,j+1}\ge x_{ij} \ge x_{i+1,j+1}$, for $1\le i \le j\le n-1$.
\end{itemize}
\end{definition}

We always depict a GT-pattern $(x_{ij})_{1\le i\le j\le n}$ by arranging the entries as follows:
\[
\begin{matrix}
x_{1n} &        & \cdots &        & \cdots &           & \cdots &           & x_{nn} \\
       & \ddots &        & \ddots &        & \revddots &        & \revddots &        \\
       &        & x_{13} &        & x_{23} &           & x_{33} &           &        \\
       &        &        & x_{12} &        & x_{22}    &        &           &        \\
       &        &        &        & x_{11} &           &        &           &
\end{matrix}.
\]
In this arrangement, the inequalities in Definition \ref{GT-patterns}
state that each entry is non-negative, and each entry not in the top
row is weakly less than its upper-left neighbor and weakly greater
than its upper-right neighbor.  We refer to the elements
$x_{1j}, \dotsc, x_{jj}$ as the $j$th row, \textit{i.e.}, the $j$th
row \emph{counted from the bottom}.  The solutions of these
inequalities define a polyhedral cone in $\R^{\frac{n(n+1)}{2}}$.  See the top of Figure
\ref{fig:tableauxbij} for an example of a GT-pattern.

\begin{definition}{\label{GT-polytopes}}
Given $\lambda,\mu \in \Z^n$, the \emph{Gelfand--Tsetlin polytope} $GT(\lambda, \mu) \subset X_n$ is the convex polytope of GT-patterns $(x_{ij})_{1\le i\le j\le n}$ satisfying the equalities
\begin{itemize}
\item $x_{in}=\lambda_i$, for $1\le i\le n$;
\item $x_{11}=\mu_1$; and  $\sum_{i=1}^j x_{ij}-\sum_{i=1}^{j-1}x_{i,j-1}=\mu_j$, for $2\le j\le n$.
\end{itemize}
In other words, $GT(\lambda, \mu)$ is the set of all GT-patterns in
$X_n$ in which the top row is $\lambda$ and the sum of the entries in
the $j$th row is $\sum_{i=1}^j \mu_j$ for $1 \leq j \leq n$. Note that when we speak of a GT-polytope $GT(\lambda, \mu)$, we assume that $\lambda$ and $\mu$ are integral.
\end{definition}

The importance of GT-polytopes stems from a classic result of
I. M. Gelfand and M. L. Tsetlin in \cite{gt}, which states that the
number of integral lattice points in the Gelfand--Tsetlin polytope
$GT(\lambda, \mu)$ equals the dimension of the weight $\mu$ subspace
of the irreducible representation of $\mathfrak{gl}_n \C$ with highest
weight $\lambda$.  These subspaces are indexed by the set
$SSYT(\lambda,\mu)$ of semi-standard Young tableaux with shape
$\lambda$ and content $\mu$ \cite{stanleyII}. It is well-known that
the elements of $SSYT(\lambda,\mu)$ are in one-to-one correspondence
with the integral GT-patterns in $GT(\lambda, \mu)$ under the
bijection exemplified in Figure \ref{fig:tableauxbij}: Given an
integral GT-pattern in $X_n$, let $\lambda^{(j)}$ be the $j$th row (so
that $\lambda^{(n)} = \lambda$).  For $1 \leq j \leq n$, place $j$'s
in each of the boxes in the skew shape $\lambda^{(j)}/\lambda^{(j-1)}$
in the Young diagram of shape $\lambda$.  (Here we put $\lambda^{(0)}
= \emptyset$ to deal with the $j=1$ case.) See \cite{stanleyII} for
details and \cite{berensteinkirillov,kirillov2} for more interesting
uses of Gelfand--Tsetlin polytopes.  Now we introduce the main
combinatorial tool for the study of vertices of the Gelfand--Tsetlin
polytopes:

\begin{figure}
\psfrag{x}[cc][cc][1][0]{$\begin{array}{ccccccccccc}
6 &   & 3 &   & 2 &   & 2 &   & 0 &   \\
  & 4 &   & 2 &   & 2 &   & 0 &   &   \\
  &   & 4 &   & 2 &   & 1 &   &   &   \\
  &   &   & 3 &   & 1 &   &   &   &   \\
  &   &   &   & 3 &   &   &   &   &
\end{array}
$}
\psfrag{1}[cl][cl][1][0]{1}
\psfrag{2}[cl][cl][1][0]{2}
\psfrag{3}[cl][cl][1][0]{3}
\psfrag{4}[cl][cl][1][0]{4}
\psfrag{5}[cl][cl][1][0]{5}
\begin{center}
\includegraphics[scale=.4]{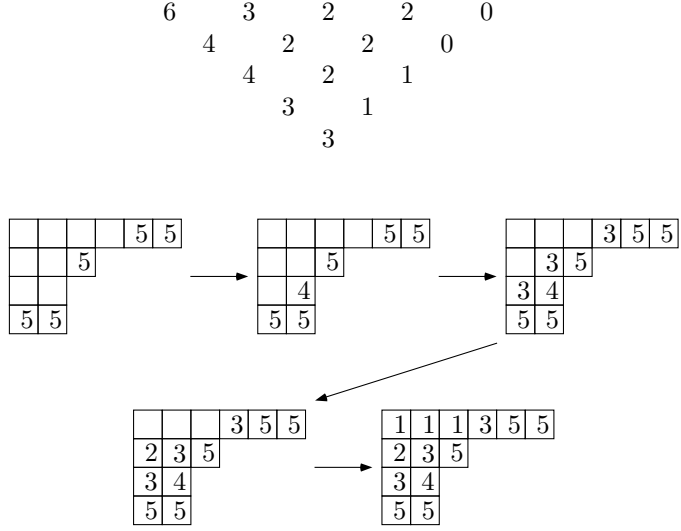}
\caption{A bijection mapping $GT(\lambda, \mu) \cap Z^{(n(n+1)/2)} \to SSYT(\lambda,\mu)$}
\label{fig:tableauxbij}
\end{center}
\end{figure}

\begin{definition}{\label{tilings}}
Given a GT-pattern $\vec{x} \in X_n$, the \emph{tiling} $\mathscr{P}$ of $\vec{x}$ is the partition of the set
\[
\{ (i,j) \in \Z^2 \colon 1 \leq i \leq j \leq n\}
\]
into subsets, called \emph{tiles}, that results from grouping together those entries in $\vec{x}$ that are equal and adjacent.  More precisely, $\mathscr{P}$ is that partition of $\{ (i,j) \in \Z^2 \colon 1 \leq i \leq j \leq n\}$ such that two pairs $(i,j), (\tilde{i},\tilde{j})$ are in the same tile iff there are sequences
\begin{gather*}
i = i_1, i_2, \ldots, i_r = \tilde{i}, \\
j = j_1, j_2, \ldots, j_r = \tilde{j}
\end{gather*}
such that $x_{i_{k+1} j_{k+1}} = x_{i_k j_k}$ and
\[
(i_{k+1}, j_{k+1}) \in \{ (i_k + 1, j_k + 1), (i_k, j_k + 1), (i_k - 1, j_k - 1), (i_k, j_k - 1) \},
\]
for each $k \in \{1, \ldots, r-1 \}$.
\end{definition}

Given a GT-pattern $\vec{x}$ with tiling $\mathscr{P}$, we associate
to $\mathscr{P}$ (or, equivalently, to $\vec{x}$) a matrix $A_\mathscr{P}$ as follows.  
Define the \emph{free tiles} $P_1, P_2, \ldots, P_s$ of $\mathscr{P}$ to be those tiles in
$\mathscr{P}$ that do not intersect the bottom or top row of
$\vec{x}$, \textit{i.e.}, those tiles that do not contain $(1,1)$ and
do not contain $(i,n)$ for $1 \leq i \leq n$.  The order in which the free tiles are indexed will not matter for our purposes, but, for concreteness, we adopt the convention
of indexing the free tiles in the order that they are initially
encountered as the entries of $\vec{x}$ are read from left to right
and bottom to top.  Define the \emph{tiling matrix} $A_{\mathscr{P}} = (a_{jk})_{2 \leq j \leq
n-1, \, 1 \leq k \leq s}$ by
\[
a_{jk} = \# \{i \colon (i,j) \in P_k \}.
\]
(Note that the index $j$ begins at 2.)  That is, $a_{jk}$ counts the number of entries in the $j$th row of $\vec{x}$ that are contained in $P_k$.

\begin{figure}
\psfrag{0}[cc][cc][1][0]{$0$}
\psfrag{1}[cc][cc][1][0]{$1$}
\psfrag{2}[cc][cc][1][0]{$2$}
\psfrag{3}[cc][cc][1][0]{$3$}
\psfrag{4}[cc][cc][1][0]{$4$}
\psfrag{5}[cc][cc][1][0]{$5$}
\psfrag{6}[cc][cc][1][0]{$6$}
\psfrag{1o2}[cc][cc][1][0]{$\frac{1}{2}$}
\psfrag{7o2}[cc][cc][1][0]{$\frac{7}{2}$}
\psfrag{9o2}[cc][cc][1][0]{$\frac{9}{2}$}
\begin{center}
\includegraphics[scale=.4]{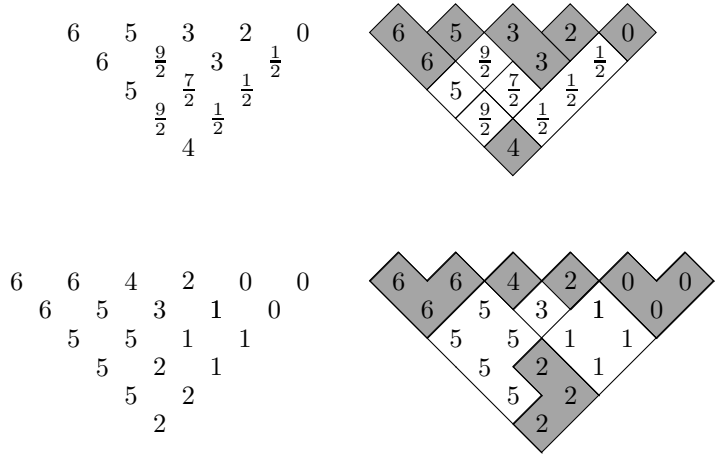}
\caption{Tilings of GT-patterns}
\label{fig:tiling}
\end{center}
\end{figure}

\begin{exm}
Two GT-patterns and their tilings are given in Figure \ref{fig:tiling}.  The unshaded tiles are the free tiles.  The associated tiling matrices are respectively
\[
\begin{bmatrix}
1 & 1 & 0 & 0 & 0 \\
0 & 1 & 1 & 1 & 0 \\
0 & 1 & 0 & 0 & 1
\end{bmatrix}
\]
and
\[
\begin{bmatrix}
1 & 0 & 0 \\
1 & 1 & 0 \\
2 & 2 & 0 \\
1 & 1 & 1 
\end{bmatrix}.
\]
\end{exm}

The motivation for introducing tilings, and the main result of this paper, is the following.

\begin{thm}{\label{facedimension}}
Suppose that $\mathscr{P}$ is the tiling of a GT-pattern $\vec{x}$.  Then the dimension of the kernel of $A_\mathscr{P}$ is equal to the dimension of the minimal (dimensional) face of the GT-polytope containing $\vec{x}$.
\end{thm}

As a corollary to this result, we get an easy-to-check criterion for a GT-pattern being a vertex of the GT-polytope containing it.

\begin{cor}{\label{vertexcriterion}}
If $\vec{x} \in GT(\lambda,\mu)$ has tiling $\mathscr{P}$ containing
$s$ free tiles, then the following conditions are equivalent:
\begin{itemize}
\item $\vec{x}$ is a vertex of $GT(\lambda,\mu)$; and
\item $A_\mathscr{P}$ has trivial kernel; \textit{i.e}, for some $s \times s$ submatrix $\tilde{A}$ of $A_\mathscr{P}$, $\det \tilde{A} \neq 0$.
\end{itemize}
\end{cor}

As an application of Theorem \ref{facedimension}, we present a
solution to a conjecture by Berenstein and Kirillov (p.101,
Conjecture~2.1 in \cite{berensteinkirillov}): all vertices of a
Gelfand--Tsetlin polytope have integer coordinates, i.e. $GT(\lambda
,\mu)$ is a convex integral polytope. This conjecture seems to have
been motivated by the fact that, for an integer parameter $l$, the
Kostka number $K_{l\lambda ,l\mu}$ is a polynomial in $l$.  This was
proved by Kirillov and Reshetikin using some fermionic formulas in
\cite{kirillov}. More recently, Billey et al.  \cite{billeyetal}
presented a more direct proof.  For completeness, we give another
proof at the end of Section \ref{proofofresults}. Derksen and Weyman
\cite{derksenweyman} have now extended this polynomiality to other Ehrhart
functions, the so-called Littlewood-Richardson coefficients. 

We must comment that it is quite natural to conjecture integrality of
the vertices of GT-polytopes, if one knows of the theory of Ehrhart
functions that count the number of lattice points inside convex
polytopes and their dilations (see Chapter 4 \cite{stanleyI}). The
Ehrhart counting functions are known to be polynomials when the
vertices are integral. As a consequence, in the following theorem we
are in fact presenting the first infinite family of non-integral
polyhedra whose Ehrhart counting functions are still
polynomials. Other low-dimensional families have been found recently
\cite{mcallisterwoods}. Finally, we must remark that R.P. Stanley
communicated to us that his student Peter Clifford noticed
non-integrality for GT-polytopes earlier (unwritten) and that King et
al. had independently noticed non-integrality for hive polytopes
(which generalize the GT-polytopes) in the case $n=5$ (see
\cite{kingetal}).  They also proved integrality of vertices for $n
\leq 4$, did a nice study of ``stretched'' Kostka coefficients, and
presented several conjectures concerning again polynomiality of
Ehrhart counting functions.

\begin{thm}{\label{counterexamplesBK}}
The Berenstein--Kirillov conjecture is true for $n \leq 4$, but
counterexamples to this conjecture exist for all values of $n \geq 5$.
More strongly, by choosing $n$ sufficiently large, we can find
GT-polytopes in which the denominators of the vertices are arbitrarily
large: For positive integer $k$, let $\lambda = (k^k, k-1, 0^k)$ and
$\mu = ((k-1)^{k+1},1^k)$.  Then a vertex of $GT(\lambda, \mu) \subset
X_{2k+1}$ contains entries with denominator $k$.
\end{thm}

\section{Proof of the Main Result and its Consequences} \label{proofofresults}

\begin{proof}[Proof of Theorem \ref{facedimension}]  Suppose that $\mathscr{P}$ is the tiling of a GT-pattern $\vec{x}$ in the GT-polytope $GT(\lambda, \mu) \subset X_n$.  Let $s$ be the number of free tiles in $\mathscr{P}$.
Let $(\epsilon^{(1)}, \dotsc, \epsilon^{(d)})$ be a basis for $\ker
A_\mathscr{P}$. Because we can scale the basis by any nonzero scalar,
we can assume that 
\[
\abs{\epsilon^{(m)}_k} < 1/2 \min \{ \abs{x_{i_1 j_1} - x_{i_2 j_2}} \colon x_{i_1 j_1} \neq x_{i_2 j_2} \}, \qquad \text{for $1 \leq m \leq d$, $1 \leq k \leq s$,} 
\]
where $\epsilon^{(m)}_k$ is the $k$th coordinate of $\epsilon^{(m)}$.

Let $H \subset X_n$ be the linear subspace of $X_n$ such that $H + \vec{x}$ is the affine span of the minimal face of $GT(\lambda, \mu)$ containing $\vec{x}$.  Define a linear map $\varphi \colon \ker A_\mathscr{P} \to X_n$ by $\varphi(\epsilon^{(m)}) = \vec{y}^{(m)}$, where
\[
y^{(m)}_{ij} = \begin{cases}
\epsilon^{(m)}_k & \text{if $(i,j)$ is in the free tile $P_k$ of $\mathscr{P}$}, \\
0 & \text{if $(i,j)$ is not in a free tile of $\mathscr{P}$}.
\end{cases}
\]
(See Example \ref{example2}) Thus, $\vec{x} + \vec{y}^{(m)}$ is the result of adding $\epsilon^{(m)}_k$ to each entry in the $k$th free tile of $\vec{x}$ for $1 \leq k \leq s$.

The claim is that $(\vec{y}^{(1)}, \dotsc, \vec{y}^{(d)})$ is a basis
for $H$.  First, since the $\epsilon^{(m)}_k$'s are sufficiently
small, $\vec{x} \pm \vec{y}^{(m)}$ is a GT-pattern.  Moreover,
$y^{(m)}_{11} = 0$, $y^{(m)}_{in}=0$ for $1 \leq i \leq n$, and each
row-sum of $\vec{y}^{(m)}$ is 0. This last fact is true because
$\epsilon^{(m)} \in \ker A_\mathscr{P}$ and the row-sum is, by
construction, the same as the dot product of $\epsilon^{(m)}$ with the
matrix $A_\mathscr{P}$. Taken together, these properties yield that $\vec{x} \pm \vec{y}^{(m)} \in
GT(\lambda, \mu)$.  That is, $\vec{x} + \vec{y}^{(m)}$ and $\vec{x} -
\vec{y}^{(m)}$ are the endpoints of a line segment contained in
$GT(\lambda,\mu)$ that contains $\vec{x}$ in its relative interior.
This establishes that $\vec{y}^{(1)}, \dotsc, \vec{y}^{(d)} \in H$.

That $\vec{y}^{(1)}, \dotsc, \vec{y}^{(d)}$ are linearly independent
clearly follows from the fact that $\epsilon^{(1)}, \dotsc,
\epsilon^{(d)}$ are linearly independent.  Thus, it remains only to
prove that $\vec{y}^{(1)}, \dotsc, \vec{y}^{(d)}$ span $H$.  Suppose that
$\vec{y} \in H$, and assume that $\vec{y}$ is scaled by a nonzero amount so that 
$\vec{x} \pm \vec{y} \in GT(\lambda, \mu)$.  We construct
an element $\epsilon$ of $\ker A_\mathscr{P}$ such that $\varphi(\epsilon) = \vec{y}$.  Note that
\begin{itemize}
\item $y_{ij} = 0$ when $(i,j)$ is not in a free tile of $\mathscr{P}$,
\item each row-sum of $\vec{y}$ is 0, and
\item if $(i_1, j_1)$ and $(i_2, j_2)$ are in the same free tile of $\mathscr{P}$, then $y_{i_1 j_1} = y_{i_2 j_2}$.
\end{itemize}
To see that this last property holds, it suffices (see Definition (\ref{tilings})) to examine the case where $y_{i_1 j_1}$ and $y_{i_2 j_2}$ are adjacent entries, \textit{i.e.} where
\[
(i_2, j_2) \in \{ (i_1 + 1, j_1 + 1), (i_1, j_1 + 1), (i_1 - 1, j_1 - 1), (i_1, j_1 - 1) \}.
\]
Since $\vec{x} \pm \vec{y}$ is a GT-pattern (see Definition (\ref{GT-patterns})), we must have either
\[
x_{i_1 j_1} + y_{i_1 j_1} \leq x_{i_2 j_2} + y_{i_2 j_2} \qquad \text{and} \qquad x_{i_1 j_1} - y_{i_1 j_1} \leq x_{i_2 j_2} - y_{i_2 j_2}
\]
or
\[
x_{i_1 j_1} + y_{i_1 j_1} \geq x_{i_2 j_2} + y_{i_2 j_2} \qquad \text{and} \qquad x_{i_1 j_1} - y_{i_1 j_1} \geq x_{i_2 j_2} - y_{i_2 j_2}.
\]
But since $(i_1, j_1)$ and $(i_2, j_2)$ are in the same tile of $\mathscr{P}$, we have $x_{i_1 j_1} = x_{i_2 j_2}$.  Thus, in either case, we can subtract the $\vec{x}$ entries from both sides, yielding $y_{i_1 j_1} = y_{i_2 j_2}$, as claimed.

For $1\leq k \leq s$ and for each $(i,j)$ in the free tile $P_k$, put $\epsilon_k = y_{ij}$.  Let $\epsilon = (\epsilon_1, \ldots, \epsilon_s)$.  Then, from the conditions on $\vec{y}$ given above, $\epsilon \in \ker A_\mathscr{P}$ and $\varphi(\epsilon) = \vec{y}$.  Hence, the coordinates of $\epsilon$ with respect to the basis $(\epsilon^{(1)}, \dotsc, \epsilon^{(d)})$ of $\ker A_\mathscr{P}$ will also be the coordinates of $\vec{y}$ with respect to $(\vec{y}^{(1)}, \dotsc, \vec{y}^{(d)})$.  In particular, $(\vec{y}^{(1)}, \dotsc, \vec{y}^{(d)})$ is a basis for $H$, as claimed.
\end{proof}

\begin{exm}{\label{example2}}
Let $\vec{x}$ be the GT-pattern
\[
\begin{matrix}
6 &   & 5 &             & 3           &             & 2           &             & 0 \\
  & 6 &   & \frac{9}{2} &             & 3           &             & \frac{1}{2} &   \\
  &   & 5 &             & \frac{7}{2} &             & \frac{1}{2} &             &   \\
  &   &   & \frac{9}{2} &             & \frac{1}{2} &             &             &   \\
  &   &   &             & 4           &             &             &             &
\end{matrix}
\]
from Figure \ref{fig:tiling}.  We explicitly apply to $\vec{x}$ the constructions in the proof of Theorem \ref{facedimension}.  This GT-pattern has tiling matrix
\[
A_\mathscr{P} = \begin{bmatrix}
1 & 1 & 0 & 0 & 0 \\
0 & 1 & 1 & 1 & 0 \\
0 & 1 & 0 & 0 & 1
\end{bmatrix}.
\]
A ``sufficiently short" basis for $\ker A_\mathscr{P}$ is
\[
(\epsilon^{(1)}, \epsilon^{(2)}) =
\left( 1/3 \begin{bmatrix}
0  \\
0  \\
-1 \\
1  \\
0
\end{bmatrix},
1/3 \begin{bmatrix}
1  \\
-1 \\
1  \\
0  \\
1
\end{bmatrix} \right).
\]
Therefore, $\vec{x}$ lies in a 2-dimensional face of
\[
GT((6,5,3,2,0),(4,1,4,5,2)).
\]
Applying the map $\varphi$ from the proof to $(\epsilon^{(1)}, \epsilon^{(2)})$ yields
\[
\vec{y}^{(1)} =
\begin{matrix}
0 &   & 0            &   & 0           &   & 0 &   & 0 \\
  & 0 &              & 0 &             & 0 &   & 0 &   \\
  &   & -\frac{1}{3} &   & \frac{1}{3} &   & 0 &   &   \\
  &   &              & 0 &             & 0 &   &   &   \\
  &   &              &   & 0           &   &   &   &
\end{matrix}
\]
and
\[
\vec{y}^{(2)} = 
\begin{matrix}
0 &   & 0           &             & 0 &              & 0            &              & 0 \\
  & 0 &             & \frac{1}{3} &   & 0            &              & -\frac{1}{3} &   \\
  &   & \frac{1}{3} &             & 0 &              & -\frac{1}{3} &              &   \\
  &   &             & \frac{1}{3} &   & -\frac{1}{3} &              &              &   \\
  &   &             &             & 0 &              &              &              &
\end{matrix}.
\]
From the proof just given, the affine subspace affinely spanned by the face containing $\vec{x}$ is affinely spanned by $\{ \vec{x}, \vec{x} + \vec{y}^{(1)}, \vec{x} + \vec{y}^{(2)} \}$.
\end{exm}

The machinery of tilings allows us easily to find nonintegral vertices of GT-polytopes by looking for a tiling with a tiling matrix satisfying certain properties given below.  Then the tiling can be ``filled" in a systematic way with the entries of a GT-pattern that is a nonintegral vertex.

\begin{lemma}{\label{nonintvertexcriterion}}
Suppose that $\mathscr{P}$ is a tiling with $s$ free tiles such that $A_{\mathscr{P}}$ has trivial kernel.  Then the following conditions are equivalent:
\begin{enumerate}
\item $\mathscr{P}$ is the tiling of a nonintegral vertex $\vec{x}$ of a GT-polytope in which $q \in \N$ is the least common multiple of the denominators of the entries in $\vec{x}$ (written in reduced form); and
\item there is an integral vector $\xi = (\xi_1, \dotsc, \xi_s)$ such that $A_{\mathscr{P}} \xi \equiv 0(\text{mod}\,q)$ and such that, for some $k \in \{1, \dotsc, s\}$, $\gcd(\xi_k,q) = 1$.
\end{enumerate}
\end{lemma}

\begin{proof}[Proof] [(1)
$\Rightarrow$ (2)] Suppose that $\vec{x}$ is a nonintegral vertex in
which $q$ is the least common multiple of the denominators of the
entries.  For each entry $x_{ij}$, $1 \leq i \leq j \leq n$, let
$p_{ij} = q x_{ij}$.  Let $P_1, \ldots, P_s$ be the free tiles of
$\mathscr{P}$, and define $\xi = (\xi_1, \ldots , \xi_s)$ by $\xi_k =
p_{ij}$ for some $(i,j) \in P_k$ (all values of $p_{ij}$ are equal
within a tile).  Since $\vec{x}$ has entries with denominator $q$, we
have that, for some $k \in \{1, \dotsc, s\}$, $\gcd(\xi_k, q) =
1$. Moreover, since each row-sum of $\vec{x}$ is an integer, we have
that, for \emph{fixed} $j \in \{1, \ldots, n\}$,
\[
q \quad \text{divides} \quad \sum_{\substack{1 \leq k \leq s \\ (i,j) \in P_k}} p_{ij} = \sum_{1 \leq k \leq s} a_{jk} \xi_k.
\]
Therefore, $A_\mathscr{P} \xi \equiv 0(\text{mod}\,q)$.

[(2) $\Rightarrow$ (1)] $\mathscr{P}$ is given to be a tiling, so some GT-pattern $\tilde{\tilde{\vec{x}}}$ with rational entries has tiling $\mathscr{P}$.  If necessary, multiply $\tilde{\tilde{\vec{x}}}$ by some integer to produce an integral GT-pattern $\tilde{\vec{x}}$ with tiling $\mathscr{P}$.  Choose $\xi = (\xi_1, \dotsc, \xi_s)$ satisfying condition (2) such that $0 \leq \xi_1, \dotsc, \xi_s < q$.  Define $\vec{y} \in X_n$ by
\[
y_{ij} = \begin{cases}
\xi_k / q & \text{if $(i,j)$ is in the free cell $P_k$ of $\mathscr{P}$,} \\
0         & \text{if $(i,j)$ is not in a free cell of $\mathscr{P}$.}
\end{cases}
\]
Then $\vec{x} = \tilde{\vec{x}} + \vec{y}$ satisfies condition (1).
\end{proof}

Now we are ready to give the details of the proof of Theorem
\ref{counterexamplesBK}. In particular, Propositions \ref{GTnleq4} and
\ref{GTdenomsunbounded} settle the Berenstein-Kirillov
conjecture. Proposition \ref{GTnleq4} has also been proven by King et
al. \cite{kingetal} using the hive model of Knutson and Tao
\cite{knutsontao2}. We
give here a ``tiling'' proof.

\begin{prop}{\label{GTnleq4}}
When $n \leq 4$, every GT-polytope in $X_n$ is integral.
\end{prop}

\begin{proof}  Note that it suffices to prove the $n=4$ case since there is a natural embedding $X_n \hookrightarrow X_{n+1}$ defined by $\vec{x} \mapsto \tilde{\vec{x}}$, where
\[
\tilde{x}_{ij} = \begin{cases} 0 & \text{if $1 \leq i = j \leq n+1$},
\\ x_{i,j-1} & \text{if $1 \leq i < j \leq n+1$}.
\end{cases}
\]

Suppose that $\vec{x} \in X_4$ is a vertex.  Then, by Corollary \ref{vertexcriterion}, the associated tiling matrix $A_\mathscr{P}$ has trivial kernel.  Therefore, $A_\mathscr{P}$ is either a $2\times1$ or a $2\times2$ matrix.  Note also that the first and last nonzero entries of each column of a tiling matrix associated with a GT-pattern must be 1.  Therefore, $A_\mathscr{P}$ is a 0/1-matrix.

If $A_\mathscr{P}$ is $2\times1$, then the only possibilities are 
\[
A_\mathscr{P} = \left[\begin{array}{c}
1 \\
0
\end{array}\right], \,
A_\mathscr{P} = \left[\begin{array}{c}
1 \\
1
\end{array}\right], \, \text{or }
A_\mathscr{P} = \left[\begin{array}{c}
0 \\
1
\end{array}\right].
\]
In each case, there exists no vector $\xi \not\equiv 0(\text{mod}\,q)$ such that $A_\mathscr{P} \xi \equiv 0(\text{mod}\,q)$ for $q > 1$, so Lemma \ref{nonintvertexcriterion} implies that the entries of $\vec{x}$ are integral.  On the other hand, if $A_\mathscr{P}$ is $2\times2$, then  $\det A_\mathscr{P} \in \{-1, 1\}$, \textit{i.e.}, $\gcd(\det A_\mathscr{P}, q) = 1$ for $q>1$.  Therefore, $A_\mathscr{P}$, considered as a linear operator on $\Z/q\Z \times \Z/q\Z$, is invertible for $q>1$, so, by Lemma \ref{nonintvertexcriterion}, $\vec{x}$ is integral.
\end{proof}

Now we show that nonintegral GT-polytopes exist in $X_n$ for each $n \geq
5$.  Moreover, by choosing $n$ sufficiently large, we can find
GT-polytopes in which the denominators of the vertices are arbitrarily
large.

\begin{prop}{\label{GTdenomsunbounded}}
For positive integer $k$, let $\lambda = (k^k, k-1, 0^k)$ and $\mu = ((k-1)^{k+1},1^k)$.  Then a vertex of $GT(\lambda, \mu) \subset X_{2k+1}$ contains entries with denominator $k$.
\end{prop}

\begin{figure}
\psfrag{k}[cc][cc][1][0]{$k$}
\psfrag{0}[cc][cc][1][0]{$0$}
\psfrag{e}[cc][cc][1][0]{$\cdots$}
\psfrag{d}[cc][cc][1][135]{$\cdots$}
\psfrag{r}[cc][cc][1][45]{$\cdots$}
\psfrag{k-}[cc][cc][1][0]{$k-1$}
\psfrag{km}[cc][cc][1][0]{$k-\frac{1}{k}$}
\begin{center}
\includegraphics[scale = .65]{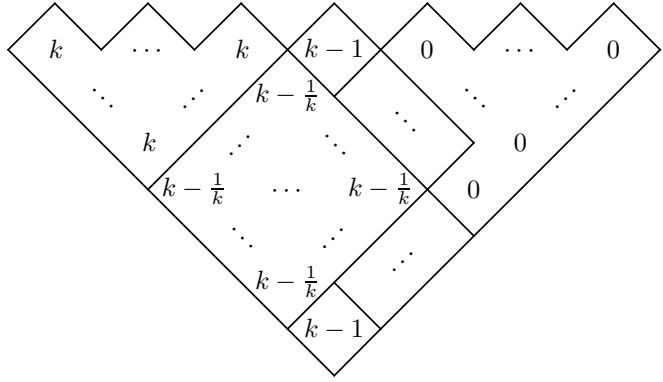}
\caption{An infinite family of counterexamples to the Berenstein-Kirillov conjecture}
\label{fig:counterexample}
\end{center}
\end{figure}

\begin{proof}
Define $\vec{x}^{(k)} \in X_{2k+1}$ by
\[
x_{ij}^{(k)} = \begin{cases}
\frac{(k-j+1)(k+1)}{k} & \text{if $1 \leq i = j \leq k+1$,}                     \\
k - \frac{1}{k}        & \text{if $1 \leq i < j \leq k+1$,}                     \\
k                      & \text{if $k+1 < j \leq 2k+1$ and $1 \leq i < j-k$,}    \\
k - \frac{1}{k}        & \text{if $k+1 < j \leq 2k+1$ and $j-k \leq i \leq k$,} \\
\frac{(j-k-1)(k-1)}{k} & \text{if $k+1 < j \leq 2k+1$ and $i = k+1$,}           \\
0                      & \text{if $k+1 < j \leq 2k+1$ and $k+1 < i \leq 2k+1$.}
\end{cases}
\]
(See Figure \ref{fig:counterexample}.)  Then $\vec{x}^{(k)} \in GT(\lambda,\mu)$.  The tiling matrix associated with $\vec{x}^{(k)}$ is
\[
A_\mathscr{P} = \begin{bmatrix}
1      & 1      & 0      & \cdots & 0      & 0      & \cdots & 0      & 0      \\
2      & 0      & 1      & \cdots & 0      & 0      & \cdots & 0      & 0      \\
\vdots & \vdots & \vdots & \ddots & \vdots & \vdots & \ddots & \vdots & \vdots \\
k - 1  & 0      & 0      & \cdots & 1      & 0      & \cdots & 0      & 0      \\
k      & 0      & 0      & \cdots & 0      & 0      & \cdots & 0      & 0      \\
k - 1  & 0      & 0      & \cdots & 0      & 1      & \cdots & 0      & 0      \\
\vdots & \vdots & \vdots & \ddots & \vdots & \vdots & \ddots & \vdots & \vdots \\
2      & 0      & 0      & \cdots & 0      & 0      & \cdots & 1      & 0      \\
1      & 0      & 0      & \cdots & 0      & 0      & \cdots & 0      & 1
\end{bmatrix}.
\]
Since $\det A_\mathscr{P} = k$, $\vec{x}^{(k)}$ is a vertex of $GT(\lambda,\mu)$.
\end{proof}

Proposition \ref{GTdenomsunbounded} explicitly constructs
counterexamples to the Berenstein-Kirillov conjecture in $X_n$ where
$n \geq 5$ is odd.  Counterexamples with even $n \geq 6$ may be
constructed from these using the embedding $X_n \hookrightarrow
X_{n+1}$ given in the proof of Theorem \ref{GTnleq4}.  Less trivial
examples with even $n$ may be constructed using other tilings.

It is interesting to note that if we fix the size of the
GT-patterns, then arbitrarily large denominators in the vertices of
GT-polytopes do not exist. To see this, we observe that Lemma
\ref{nonintvertexcriterion} says that if $\vec{x}$ is a nonintegral
vertex in which $q$ appears as a denominator, then the tiling matrix
$A_{\mathscr{P}}$ has trivial kernel as a linear operator $\R^s \to
\R^{n-2}$ (since $\vec{x}$ is a vertex), but $A_{\mathscr{P}}$ has
\emph{nontrivial} kernel when considered as an operator $(\Z/q\Z)^s
\to (\Z/q\Z)^{n-2}$. Moreover, this nontrivial kernel contains a
vector in which one of the coordinates is a unit in $\Z/q\Z$.  This
last condition implies that each $s \times s$ submatrix of
$A_{\mathscr{P}}$ has determinant equal to 0 modulo $q$.

\begin{prop}{\label{denombound}}
For fixed $n$, the numbers that may appear as denominators of entries
in vertices of GT-polytopes in $X_n$ are smaller than $(n-1)^{ {n+1
\choose 2} -n-1}$.
\end{prop}

\begin{proof}

Fix $n \in \N$.  Since only finitely many partitions of $\{ (i,j) \in
\Z^2 \colon 1 \leq i \leq j \leq n\}$ exist, there is an upper bound
on the set
\[
\left\lbrace \abs{m} \colon
\begin{array}{c}
\text{ $m$ is the determinant of a square row submatrix          } \\
\text{ of the tiling matrix of some GT-pattern $\vec{x} \in X_n$ }
\end{array} \right\rbrace.
\]

Let $N$ be an upper bound on this set.  The claim is that no
GT-polytope in $X_n$ has a vertex with denominators greater than $N$.

Let $q > N$ be given.  Suppose that $\vec{x} \in X_n$ is a vertex.  Let $s$ be the number of free tiles in $\vec{x}$, and let $A_\mathscr{P}$ be the tiling matrix of $\vec{x}$.  Then no $s \times s$ submatrix of $A_\mathscr{P}$ has determinant greater than or equal to $q$.  Moreover, by Corollary \ref{vertexcriterion}, some $s \times s$ submatrix of $A_\mathscr{P}$ has nonzero determinant.  Therefore, this $s \times s$ submatrix has determinant not equal to 0 modulo $q$.  However, in the remarks preceding this proposition, we noted that if $\vec{x}$ is a vertex in which $q$ is a denominator of one of the entries, then \emph{every} $s \times s$ submatrix has determinant equal to 0 modulo $q$.  This proves that N is a bound as claimed.

Our second claim is that $N$ is no more than $(n-1)^{{n+1
\choose 2} -n-1}$. We need to bound the $s \times s$ subdeterminants.
In a tiling matrix all entries are nonnegative integers, so we 
know that the spectral radius of eigenvalues possible for such a
matrix is bounded above by the maximum possible sum of entries
along a row \cite{hornjohnson}. A row sum of the tiling matrix
cannot be more than the number of elements in a middle row of 
a GT-pattern, and this is $n-1$ for the longest row. On the other
hand, how big can $s$ be? It is no bigger than the number
of entries in the middle of the GT-pattern, which equals 
${n+1 \choose 2} -n-1$. Since the determinant is less than
the product of the norms of the eigenvalues, we get the desired bound.
\end{proof}

\noindent To conclude this paper we present another proof of the following
result:

\begin{prop} Given the GT-polytope $GT(\lambda, \mu) \subset X_n$, 
the Ehrhart counting function 
$f(m)=\#\left( GT(m\lambda, m\mu) \cap {\mathbb Z}^{n+1 \choose 2} \right)$ 
is a univariate polynomial.
\end{prop}

\begin{proof} 
It is well-known, from Ehrhart's fundamental work, that $f(m)$ must be
a \emph{quasipolynomial}. This means that there exist an integer $M$
and polynomials $g_0,g_1,\dots,g_{M-1}$ such that $f(m) = g_i(m)$ if $m
\equiv i \ \mod{M}$ (see details in Chapter 4 of \cite{stanleyI}).
So it is then enough to prove that, for some large enough value of $m$,
a single polynomial interpolates all values from then on, because 
then the $g_i$'s are forced to coincide infinitely many times, which
proves that they are the same polynomial.

We use the algebraic meaning of $f(m)$ as the multiplicity of weight
$n\mu$ in the irreducible representation $V_{n\lambda}$ of
$\mathfrak{gl}_n \C$. The well-known Kostant's multiplicity formula
(see page 421 of \cite{fultonharris}) gives that
$$
f(m)=\sum_{\sigma \in S_n} (-1)^{l(\sigma)}
K(\sigma(m\lambda+\delta)-m\mu-\delta), \qquad (*)
$$
where $K(b)$ is Kostant's partition function for the root system
$A_n$, $l(\sigma)$ denotes the number of inversions of $\sigma$, and
$\delta$ is one-half of the sum of positive roots in $A_n$.  

Kostant's partition function is what combinatorialists call a
\emph{vector partition function}\cite{sturmfels}.  More precisely,
$K(b)$ is equal to the number of nonnegative integral solutions $x$ of
a linear system $Ax=b$. The columns of $A$ are exactly the positive
roots of the system $A_n$. Because the matrix $A$ is unimodular
\cite{schrijver}, the counting function $K(b)$ is a multivariate
piecewise polynomial function.  The regions where $K(b)$ is a
polynomial are convex polyhedral cones called \emph{chambers}
\cite{sturmfels}. The chamber that contains $b$ determines the
polynomial value of $K(b)$; in fact it is the vector direction of $b$,
not its norm, that determines the polynomial formula to be used.

In formula $(*)$ the right-hand side vector for Kostant's partition
function is $b=\sigma(m\lambda+\delta)-(m\mu+\delta).$ As $m$ grows,
we might be moving from one chamber to another. Our claim is that, from
some value of $m$ on, the vectors $\sigma(m\lambda+\delta)-(m\mu+\delta)$ are inside the same
chamber.  To see this, note that in the expression $(*)$, $\mu,\lambda$, and
$\delta$ are constant vectors. For a given permutation $\sigma$, the
vector direction $\sigma(m\lambda+\delta)$ is closer and closer to that of
$\sigma(\lambda)$ when $m$ grows in
value. Similarly, the vector direction of $m\mu+\delta$ approaches that of $\mu$ when $m$ grows. Thus, the direction of $\sigma(m\lambda+\delta)-(m\mu+\delta)$ will approach the direction of
$b'=\lambda+\mu$ along a straight line. For sufficiently large $m$,
the 
vectors $\sigma(m\lambda+\delta)-(m\mu+\delta)$ are contained in the same chamber, then a single polynomial gives the value of
$K(b')$. 

Finally, we have that, for all values of $m$ greater than some $M$,
the formula $(*)$ represents an alternating sum of polynomials in the
variable $m$.  Therefore $f(m)$ is a polynomial for all $m$ greater
than $M$, exactly as we wished to prove.
\end{proof}


\begin{thebibliography}{99}

\bibitem{berensteinkirillov} A.~Berenstein and A.N. Kirillov, 
Groups generated by involutions, Gelfand-Tsetlin patterns, 
and combinatorics of Young tableaux, Algebra i Analiz 7 (1995),
no. 1, 92--152 (Russian). Translation in St. Petersburg Math. J. 7
(1996), no. 1, 77-127.

\bibitem{berensteinzelevinsky} A. D. Berenstein and A. V. Zelevinsky, Tensor product multiplicities and convex polytopes in partition space, {\em J. Geom. Phys.} {\bf 5} (1988), no. 3, 453-472.

\bibitem{billeyetal} S. Billey, V. Guillemin, and E. Rassart, A vector
partition function for the multiplicities of $\mathfrak{sl}_k {\mathbb
C}$, preprint 2003.

\bibitem{derksenweyman} H. Derksen and J. Weyman, On the Littlewood-Richardson polynomials,
\emph{J. of Algebra}. \textbf{255} (2002),  no. 2, 247-257.


\bibitem{kingetal} R.C. King, C. Tollu, and F. Toumazet,
Stretched Littlewood-Richardson and Kostka coefficients, to appear in
\emph{CRM Proceedings and Lecture Notes}, Vol 34, 2003.

\bibitem{kirillov2} A. N. Kirillov, Ubiquity of Kostka polynomials, {\em Physics and combinatorics 1999 (Nagoya)},  85-200, {\em World Sci. Publishing, River Edge, NJ}, 2001.

\bibitem{kirillov} A. N. Kirillov and N. Y. Reshetikhin, The Bethe Ansatz and the Combinatorics of Young Tableaux, {\em J. Soviet Math} {\bf 41} (1988), 925-955.

\bibitem{knutsontao} A. Knutson and T. Tao, The honeycomb model of $GL_n(\C)$ tensor products I: Proof of the Saturation Conjecture, {\em J. Amer. Math. Soc.} {\bf 12} (1999), 1055-1090.

\bibitem{knutsontao2} A. Knutson, T. Tao, and C. Woodward, The
honeycomb model of $GL_n(\C)$ tensor products II: Puzzles determine
facets of the Littlewood-Richardson cone, to appear in {\em
J. Amer. Math. Soc}, available at arXiv:math.CO/0107011.

\bibitem{fultonharris} W. Fulton and J. Harris, Representation Theory:
A first course, Graduate texts in Mathematics vol 129, Springer, New
York, 1991.

\bibitem{hornjohnson} R. Horn and C. R. Johnson, Matrix Analysis,
Cambridge Univ. Press, Cambridge, 1985.

\bibitem{gt} I. M. Gelfand and M. L. Tsetlin, Finite-dimensional representations of the group of unimodular matrices, \emph{Doklady Akad. Nauk SSSR} (N.S.) {\bf 71} (1950), 825-828),

\bibitem{mcallisterwoods} T. B. McAllister and K. Woods, manuscript in
preparation, 2003.

\bibitem{nakashimazelevinsky} T. Nakashima and A. Zelevinsky, Polyhedral realizations of crystal bases for quantized Kac-Moody algebras, {\em Adv. Math.} {\bf 131} (1997), no. 1, 253-278.

\bibitem{etienne}E. Rassart, A polynomiality property for
Littlewood-Richardson coefficients, manuscript 2003. 


\bibitem{schrijver} A. Schrijver, Theory of Linear and Integer Programming, John Wiley \& Sons Ltd., 1986.

\bibitem{stanleyI} R.P Stanley, Enumerative Combinatorics, Volume I, Cambridge University Press, 1997.

\bibitem{stanleyII} R.P Stanley, Enumerative Combinatorics,
Volume II, Cambridge University Press, 1999.

\bibitem{sturmfels} B. Sturmfels, On Vector Partition functions,
\emph{J. of Combinatorial Theory} series A, 72, 302--309, 1995.

\bibitem{szenesvergne} A. Szenes and M. Vergne, Residue formulae for
vector partitions and Euler-Maclaurin sums. available at
arXiv.math.CO:0202253, 2002.
\end{thebibliography}
\end{document}